\documentclass[12pt]{article}
\usepackage{amsfonts}
\usepackage{amssymb}
\usepackage{latexsym}
\usepackage{epsf,graphicx}

\newtheorem{remar}{Remark}
\newtheorem{examp}{Example}

\newcommand{\rema}{\begin{remar}\rm}
\newcommand{\erema}{$\blacktriangleright$\end{remar}\vspace{2mm}}

\newcommand{\exa}{\begin{examp}\rm}
\newcommand{\eexa}{$\blacktriangleright$\end{examp}\vspace{2mm}}

\newcommand {\Hol}{\mathop{\rm Hol}\nolimits}
\newcommand {\tr}{\mathop{\rm tr}\nolimits}

\renewcommand {\Re}{\mathop{\rm Re}\nolimits}
\newcommand{\pl}{\partial}
\newcommand{\pr}{\noindent{\bf Proof.}\quad }
\newcommand{\epr}{\ $\blacksquare$ \vspace{2mm} }

\newcommand{\be} {\begin{eqnarray}}
\newcommand{\ee} {\end{eqnarray}}
\newcommand{\bep} {\begin{eqnarray*}}
\newcommand{\eep} {\end{eqnarray*}}

\newfont{\bbb}{msbm10 at 12pt}
\def\Bbb#1{\hbox{\bbb #1}}

\newcommand{\C}{{\Bbb C}}

\newcommand{\B}{{\Bbb B}}

\newtheorem{defin}{Definition}[section]
\newtheorem{theorem}{Theorem}[section]
\newtheorem{corol}{Corollary}[section]
\newtheorem{propo}{Proposition}[section]

\newtheorem{lemma}{Lemma}[section]

\title{Extension operators via semigroups}
\author{Mark Elin
\\ {\small Department of Mathematics}
\\ {\small ORT  Braude College}
\\ {\small P.O. Box 78, 21982 Karmiel, Israel}
\\ {\small e-mail: mark$\_$elin@braude.ac.il}}

\begin{document}

\maketitle

\begin{abstract}
The Roper–-Suffridge extension operator and its modifications are
powerful tools to construct biholomorphic mappings with special
geometric properties.

The first purpose of this paper is to analyze common properties of
different extension operators and to define an extension operator
for biholomorphic mappings on the open unit ball of an arbitrary
complex Banach space. The second purpose is to study extension
operators for starlike, spirallike and convex in one direction
mappings. In particular, we show that the extension of each
spirallike mapping is $A$-spirallike for a variety of linear
operators $A$.

Our approach is based on a connection of special classes of
biholomorphic mappings defined on the open unit ball of a complex
Banach space with semigroups acting on this ball.
%
%
%
%
\end{abstract}

\section{Introduction}\label{sect-intr}
\setcounter{equation}{0}

One of the main purposes of the classical Geometric Function
Theory is the study of various classes of univalent and
multivalent mappings. Convex, starlike and spirallike functions on
the open unit disk $\Delta\in\C$ have been the objects of
intensive study for over a century. A reader can be referred to
the book of Goodman \cite{GAW}. The study of different classes of
biholomorphic mappings in multidimensional settings began later.
In fact, the first survey appeared in 1977 (see \cite{STJ-77}).
Recent developments in this area are reflected in \cite{GS-99,
G-K-book, E-R-S-2004} and \cite{RS-SD-book}. However, numerous
well-known tools for the construction of mappings with special
geometric properties on the open unit disk
$\Delta:=\left\{z\in\mathbb{C}:\ |z|<1\right\}$ have no
generalization for the multidimensional case. For example, until
recently only a few concrete examples of convex, starlike and
spirallike mappings in the open unit ball in $\C^n$ were known.

In 1995, Roper and Suffridge \cite{Ro-Su} introduced an extension
operator, which provides a variety of required examples. Given a
univalent function $f\in\Hol(\Delta,\C)$ normalized by
$f(0)=f'(0)-1=0$, they considered the mapping
$\Phi[f]:\mathbb{B}^n\mapsto\mathbb{C}^n$ defined as follows:
\begin{equation}\label{RS}
\Phi[f](z_1,x)=\left(f(z_1),\sqrt{f'(z_1)}\,x \right),
\end{equation}
where $x=\left(z_2,\ldots,z_n\right)$. The Roper--Suffridge
extension operator has remarkable properties. In particular:

\begin{itemize}
\item if $f$ is a normalized convex function on $\Delta$, then
$\Phi[f]$ is a normalized convex mapping on $\B^n$, see
\cite{Ro-Su};

\item if $f$ is a normalized starlike function on $\Delta$, then
$\Phi[f]$ is a normalized starlike mapping on $\B^n$, see
\cite{G-K-2000};

\item if $f$ is a normalized $\mu$-spirallike function on
$\Delta$, then $\Phi[f]$ is a normalized $\mu I$-spirallike
mapping on $\B^n$, see, for example, \cite{G-K-K-2000, L-L-2005}.

\item if $f$ is a normalized Bloch function on $\Delta$, then
$\Phi[f]$ is a normalized Bloch mapping on $\B^n$,
see\cite{G-K-2000}.
\end{itemize}

Several authors have discussed this operator and its
generalizations. In particular, the operator
\begin{equation}\label{G-K-K}
\Phi_\alpha[f](z_1,x)=\left(f(z_1),\left(f'(z_1)\right)^\alpha\,x
\right),
\end{equation}
where $\alpha\in[0,\frac12]$, was introduced in \cite{G-K-K-2000}.

For a locally biholomorphic mapping $f$ defined on the unit ball
of $\mathbb{C}^n$, Pfaltzgraff and Suffridge constructed in
\cite{P-S} an extension operator as follows:
\begin{equation}\label{PS}
\widehat{\Phi}_n[f](z,x)=\left(f(z),\left(J_f(z)\right)^{\frac1{n+1}}\,x
\right),
\end{equation}
where $z\in\C^n,\ x\in\mathbb{C},\ \|z\|^2+|x|^2<1,$ and $J_f(z)$
is the complex Jacobian of the mapping $f$ at the point $z$. It
was shown in \cite{G-K-P-2007} that this operator preserves
starlikeness.

Another extension operator was introduced in \cite{G-K-2000} for
locally biholomorphic functions $f\in\Hol(\Delta,\C)$ by
\begin{equation}\label{G-K}
\widetilde{\Phi}_\beta[f](z_1,x)=\left(f(z_1),\left(\frac{f(z_1)}{z_1}\right)^\beta\,x
\right),
\end{equation}
where $\beta\in[0,1]$. These extension operators and their
combinations (with multiplier
$\left(f'(z)\right)^{\alpha_j}\left(\frac{f(z)}z\right)^{\beta_j}$
in $j$-th coordinate) in the space $\C^n$ equipped with different
concrete norms have been considered in numerous papers. Detailed
references can be found in \cite{F-L}.

Note that, as we updated, all extension operators were studied for
functions $f$ satisfying the standard normalization $f(0)=0$ and
$f'(0)=1$ (or $J_f(0)=\rm id$, respectively).

The first purpose of this paper is to analyze common properties of
different extension operators and to define an extension operator
for biholomorphic mappings on the open unit ball of an arbitrary
complex Banach space.

The second purpose is to study extension operators for mappings
starlike or spirallike with respect to an arbitrary interior or a
boundary point (see definitions in Section~\ref{sect-nota}).
Although the case of spirallikeness with respect to an interior
point can often be reduced to a standard one ($f(0)=0$), extension
operators for starlike and spirallike mappings with respect to a
boundary point have not been considered at all. The following
effect is new even for the case of $f(0)=0$: we show that the
extension of each spirallike function is $A$-spirallike for a
variety of linear operators $A$.

Our approach is based on several simple but effective
observations:

(1) All extension operators mentioned above have the form:
\[
f(x)\mapsto \left(f(x),\Gamma(f,x)y \right)
\]
with a certain linear operator $\Gamma$ depending on a mapping $f$
and a point $x$. So, we have to understand which properties of
$\Gamma$ enable us to use it to construct an extension operator.
We will say that operators having such properties are appropriate.

(2) A biholomorphic mapping is $A$-spirallike if and only if its
image is $S$-invariant, where $S=\left\{e^{-tA}\right\}_{t\ge0}\ $
is the semigroup of linear transformations. Similar relations
between biholomorphic mappings and special semigroups also exist
for other classes of biholomorphic mappings. Therefore, we must
study extension operators for one-parameter continuous semigroups.

(3) Extension operators for a semigroup of biholomorphic
self-mappings of the open unit ball and for a corresponding class
of biholomorphic mappings do not necessarily coincide.

\section{Preliminary notions}\label{sect-nota}
\setcounter{equation}{0}

In this section we present some notions of nonlinear analysis and
geometric function theory which will be useful subsequently. A
reader may be referred to as the book \cite{RS-SD-book}.

Let $X$ be a complex Banach space with the norm $\|\cdot\|$.
Denote by $\Hol(D,E)$ the set of all holomorphic mappings on a
domain $D\subset X$ which map $D$ into a set $E\subset X$ and a
set $\Hol(D):=\Hol(D,D)$.

We start with the notion of a one-parameter continuous semigroup.

\begin{defin}\label{def-sg}
Let $D$ be a domain in a complex Banach space $X$. A family
$S=\{F_t\}_{t\ge0}\subset\Hol(D)$ of holomorphic self-mappings of
$D$ is said to be a one-parameter continuous semigroup (in short,
semigroup) on $D$ if \be\label{semi_prop1} F_{t+s}=F_t\circ
F_s,\quad t,s\geq 0, \ee and for all $x\in D$,
\begin{eqnarray}\label{semi_prop2}        
\lim_{t\rightarrow 0^+}F_t(x) =x.
\end{eqnarray}
\end{defin}

For example, if $D$ is the unit ball of $X$ and $A\in L(X)$ is an
accretive operator, then the family
$\left\{e^{-tA}\right\}_{t\ge0}$ forms a semigroup of proper
contractions of $D$. Moreover, each uniformly continuous semigroup
of bounded linear operators can be represented by this form.

\begin{defin} \label{def-gen}
A semigroup $S=\{F_t\}_{t\ge0}$ on $D$ is said to be generated if
for each $x\in D$, there exists the strong limit
\begin{eqnarray}
f(x):=\lim_{t\rightarrow 0^{+}}\,\frac 1t\biggl(x-F_t(x)\biggr).
\end{eqnarray}

In this case the mapping $f:\ D\mapsto X$ is called the
(infinitesimal) generator of the semigroup $S$.
\end{defin}

It was established in \cite{RS-SD1} that {\sl a semigroup $S$ of
holomorphic self-mappings of $D$ is differentiable with respect to
the parameter $t\geq 0$ (hence, generated by a holomorphic
mapping) if and only if it is locally uniformly continuous on
$D$.}

The following notion connects semigroups on biholomorphically
equivalent domains.

\begin{defin}\label{def-conj}
Let $\{F_t\}_{t\ge0}$ and $\{\Psi_t\}_{t\ge0}$ be semigroups on
domains $D$ and $\Omega$ of a complex Banach space, respectively.
We say that the semigroups are conjugate if there is a
biholomorphic mapping $h:D\mapsto\Omega$ such that
\[
h\circ F_t=\Psi_t\circ h.
\]
The mapping $h$ in this relation is called the intertwining map
for the semigroups.
\end{defin}

An important class of mappings which serve intertwining maps with
semigroups of linear transformations is  the class of spirallike
mappings.

\begin{defin}[see \cite{E-R-S-2004, RS-SD-book}, cf., \cite{STJ-77, GS-99, G-K-book}]\label{def-spiral}
Let $h$ be a biholomorphic mapping defined on a domain $D$ of a
Banach space $X$. The mapping $h$ is said to be spirallike if
there is a bounded linear operator $A$ such that the function
$\Re\lambda$ is bounded away from zero on the spectrum of $A$ and
such that for each point $w\in h(D)$ and each $t\ge0$, the point
$e^{-tA}w$ also belongs to $h(D)$. In this case $h$ is called
$A$-spirallike. If $A$ can be chosen to be the identity mapping,
that is, $e^{-t}w\in h(D)$ for all $w\in h(D)$ and all $t\ge0$,
then $h$ is called starlike.
\end{defin}

In other words, a biholomorphic mapping $h\in\Hol(D,X)$ is
$A$-spirallike if and only if it intertwines some semigroup on $D$
with the semigroup $\{e^{-At}\}_{t\ge0}$.

\rema\label{rem7} For mappings defined on the direct product
$Z=X\times Y$ of two Banach spaces $X$ and $Y$, it is relevant to
consider a block-matrix $A=\left(
\begin{array}{cc}
  A_{11} & A_{12} \\
  A_{21} & A_{22} \\
\end{array}
\right)$ with operators $A_{11}\in L(X),\, A_{12}\in L(Y,X),\,
A_{21}\in L(X,Y)$ and $A_{22}\in L(Y)$ satisfying certain
conditions. In such situation, the notion of ``$\left(
\begin{array}{cc}
  A_{11} & A_{12} \\
  A_{21} & A_{22} \\
\end{array}
\right)$-spirallikeness'' should be understood in the same sense
of Definition~\ref{def-spiral}. \erema

It follows by this definition that if $h$ is an $A$-spirallike
mapping then $0\in\overline{h(D)}$.

\begin{itemize}

\item If $0\in h(D)$, then there is a unique point $\tau\in D$
such that $h(\tau)=0$, and we say that $h$ is spirallike
(starlike) with respect to an interior point.

\item Otherwise, if $0\in\partial h(D)$, we say that $h$ is
spirallike (starlike) with respect to a boundary point.
\end{itemize}

In the one-dimensional case, the class of spirallike functions
with respect to a boundary point was introduced in \cite{A-E-S}
(see also references therein). It turns out that for each function
$h$ of this class there is a point $\tau,\ |\tau|=1,$ such that
$\lim\limits_{r\to1^-}h(r\tau)=0.$ The same conclusion also holds
in many multidimensional situations. In fact, the validity of such
claim depends on the validity of an analog of Lindel\"of's
principle (see, for instance, \cite{E-S3}).

Another class of mappings closely connected with dynamical systems
consists of mappings convex in one direction. These mappings
intertwine some semigroups on a given domain $D$ with semigroups
of shifts. More precisely:

\begin{defin}\label{def-conv-one}
Let $h$ be a biholomorphic mapping defined on a domain $D$ of a
Banach space $X$, and let $\tau\in X,\ \|\tau\|=1$. The mapping
$h$ is called convex in the direction $\tau$ if for each point
$w\in h(D)$ and each $t\ge0$, the point $w+t\tau$ also belongs to
$h(D)$.
\end{defin}

In the one-dimensional case, functions convex in one direction
have been studied by many authors starting from the classical work
of M.~S.~Robertson (see, for examples, \cite{GAW}). Recently, the
interest in these functions and their geometric properties has
received an impetus because of their connection with the semigroup
theory (see \cite{E-K-R-S} and references therein).

Note also that the semigroups $\left\{e^{-At} \right\}_{t\ge0}$
and $\left\{\cdot+t\tau \right\}_{t\ge0}$ which appear in
Definitions~\ref{def-spiral} and~\ref{def-conv-one} are particular
cases of the general semigroup of affine mappings
$\left\{e^{-At}\cdot +\lambda\int_0^te^{-As}\tau ds
\right\}_{t\ge0}$, where $A\in L(X),\ \lambda\ge0$ and $\tau\in
X,\ \|\tau\|=1$.

\section{Appropriate operator-valued mappings}\label{sect-approp}
\setcounter{equation}{0}

Let $X$ and $Y$ be two complex Banach spaces endowed with the
norms $\|\cdot\|_X$ and $\|\cdot\|_Y$, respectively, and let
$\mathbb{D}_1$ and $\mathbb{D}_2$ be the open unit balls in these
spaces. On the the space $Z=X\times Y$ we wish to define a norm
depending on $\|\cdot\|_X$ and $\|\cdot\|_Y$ only. Such a norm may
be defined as follows. Let $p:[0,1]\mapsto[0,1]$ be a continuous
function which satisfies the conditions:

\begin{itemize}
\item[(a)] $p(0)=1$, $p(1)=0$;

\item[(b)] $p$ is a strongly decreasing function;

\item[(c)] $p$ is convex up: $p\left(\frac{s_1+s_2}2\right)\ge
\frac12\left(p(s_1)+ p(s_2) \right)$ for all $s_1,s_2\in[0,1]$.
\end{itemize}
Then the set
\[
\mathbb{D}:=\left\{(x,y)\in\mathbb{D}_1\times \mathbb{D}_2\subset
Z: \|y\|_Y<p\left(\|x\|_X\right)\right\}
\]
is 
the open unit ball in $Z$ with respect to some norm $\|\cdot\|$.
Actually, this norm is the Minkowski functional of the set
$\mathbb{D}$. Under our assumption, $\|(x,y)\|$ is the unique
solution $\lambda\ge\|x\|_X$ of the equation $\|y\|_Y =\lambda p
\left(\frac{\|x\|_X}\lambda\right)$. Obviously, $Z$ equipped with
this norm $\|\cdot\|$ is a complex Banach space.

\vspace{3mm}

In our study of extension operators we need the notion of
appropriate operator-valued mappings. We define this in several
steps. First, we deal with self-mappings of $\mathbb{D}_1$.

\begin{defin}\label{def-approp}
Let $\widehat{\mathcal K}$ be a subset of $\Hol(\mathbb{D}_1)$
consisting of biholomorphic mappings and closed with respect to
composition, and let $\widehat{\Gamma}:\widehat{\mathcal
K}\times\mathbb{D}_1\mapsto L(Y)$ be a mapping continuous on
$\widehat{\mathcal K}$ and holomorphic on $\mathbb{D}_1$. We say
that $\widehat{\Gamma}$ is appropriate if it satisfies the
following properties:

\begin{itemize}
\item[(i)] the identity mapping ${\rm id}_X$ of the space $X$
belongs to $\widehat{\mathcal K}$, and $\,\widehat{\Gamma}({\rm
id}_X,x)={\rm id}_Y$, the identity mapping of the space $Y$;

\item[(ii)] $\widehat{\Gamma}$ satisfies the chain rule in the
sense that
$\widehat{\Gamma}(f,g(x))\widehat{\Gamma}(g,x)=\widehat{\Gamma}(f\circ
g,x)$ for all $f,g\in\widehat{\mathcal K}$ and $x\in\mathbb{D}_1$;

\item[(iii)] for each $f\in\widehat{\mathcal K}$ and
$x\in\mathbb{D}_1$, the operator $\widehat{\Gamma}(f,x)$ is
invertible;

\item[(iv)]
$\left\|\widehat{\Gamma}(f,x)\right\|_{L(Y)}\le\displaystyle\frac{p
\left(\|f(x)\|_X\right)}{p\left(\|x\|_X\right)}\,$ for all $f\in
\widehat{\mathcal K}$ and $x\in\mathbb{D}_1$.
\end{itemize}
\end{defin}

In the following examples we set $p(s)=(1-s^q)^{1/r}$, where $q,\
r\ge1$. Thus, the unit ball in the space $Z=X\times Y$ is defined
by
\[
\mathbb{D}=\left\{(x,y):\ \|x\|_X^q+\|y\|_Y^r<1 \right\}.
\]

\exa\label{ex1} Let $X=\C^n$ be the Euclidean $n$-dimensional
complex space. We consider the scalar operator
$\widehat{\Gamma}(f,x):=\left(J_f(x)\right)^\alpha{\rm id}_Y,\
{\alpha>}0$. To verify whether this operator is appropriate, first
we choose a branch of the power $\left(J_f(x)\right)^\alpha$ such
that condition (i) of Definition~\ref{def-approp} holds.
Furthermore, we denote by $\widehat{\mathcal K}$ a set consisting
of biholomorphic self-mappings of $\mathbb{D}_1$. In particular,
we can choose $\widehat{\mathcal K}=\widehat{\mathcal K^\tau}$,
the subset of $\Hol(\mathbb{D}_1)$ consisting of all biholomorphic
self-mappings of $\mathbb{D}_1$ with a fixed point $\tau\in
\overline{\mathbb{D}_1}$.

Conditions (ii) and (iii) obviously are satisfied. In addition,
$$\displaystyle
\left\|\widehat{\Gamma}(f,x)\right\|_{L(Y)}=|J_f(x)|^\alpha\le\left(
\frac{1-\|f(x)\|_X^2}{1-\|x\|_X^2} \right)^{\frac{(n+1)\alpha}2}$$
(see \cite[Lemma 1.1]{G-K-P-2007}). Therefore, condition~(iv) will
follow by the inequality
\begin{equation}\label{ineq}
\left(\frac{1-\|f(x)\|_X^2}{1-\|x\|_X^2}\right)^{\frac{(n+1)\alpha}2}
\le \frac{(1-\|f(x)\|_X^q)^{1/r}}{(1-\|x\|_X^q)^{1/r}}\,,
\end{equation}
which obviously holds for $\displaystyle \alpha=\frac2{r(n+1)}\,$
and $q=2$. To proceed, we rewrite (\ref{ineq}) as
\[
\frac{\left(1-\|f(x)\|_X^2\right)^{\frac{(n+1)\alpha}2}}{(1-\|f(x)\|_X^q)^{1/r}}
\le
\frac{\left(1-\|x\|_X^2\right)^{\frac{(n+1)\alpha}2}}{(1-\|x\|_X^q)^{1/r}}\,.
\]
Now, if all mappings in $\widehat{\mathcal K}$ satisfy $f(0)=0$,
then ${\|f(x)\|_X\le \|x\|_X}$. Taking into account that the
function
$\displaystyle\frac{\left(1-t^2\right)^{\frac{(n+1)\alpha}2}}{(1-t^q)^{1/r}}$
is increasing in $t\in(0,1)$ for $q\le2$ and
$\displaystyle\alpha\le\frac2{r(n+1)}\,$, we conclude that in this
situation inequality (\ref{ineq}) (hence, condition (iv)) holds.
\eexa

In the next example $X=\mathbb{C}$, the complex plane, and
$\mathbb{D}_1=\Delta,$ the open unit disk in $\mathbb{C}$.



\exa\label{ex2} Consider the scalar operator
$\displaystyle\widehat{\Gamma}(f,x)=\left(\frac{f(x)}{x}\right)^\beta{\rm
id}_Y$,\ ${\beta>0}$. Namely, we set $\widehat{\mathcal K}$ to be
the set of all univalent self-mappings of $\Delta$ with $f(0)=0$.
Similar to the above example, we choose a branch of the power
$\left(\frac{f(x)}{x}\right)^\beta$ such that conditions
(i)--(iii) of Definition~\ref{def-approp} hold. Condition (iv)
follows from the Schwarz Lemma: the inequality $|f(x)|\le |x|$
implies that
\[
\left\|\widehat{\Gamma}(f,x)\right\|_{L(Y)}=\left|\frac{f(x)}{x}\right|^\beta\le1\le
\left(\frac{1-|f(x)|^q}{1-|x|^q}\right)^{1/r}.
\]

As above, it is easy to modify this example for functions of the
set $\widehat{\mathcal K^\tau}$ for any $\tau\in\Delta$. \eexa

\vspace{-6mm}

\rema\label{rem3} In the introduction we mentioned papers where
combinations of extension operators (\ref{G-K-K}) and (\ref{G-K})
were studied. Obviously, such combinations are included in our
scheme; namely, we can consider non-scalar operators based on
Examples~\ref{ex1} and~\ref{ex2} above. \erema

In the next example, $X$ is a complex Hilbert space with inner
product $\langle\cdot ,\cdot\rangle$ and induced norm
$\|\cdot\|_X$, and $\tau\in\pl\mathbb{D}_1\subset X$. Also, we set
$p(s)=(1-s^2)^ {1/r}$.

\exa\label{ex3}  Consider the scalar operator
$\displaystyle\widehat{\Gamma}(f,x)=\left(\frac{1-\langle f(x),
\tau\rangle}{1-\langle x,\tau\rangle}\right)^{2/r}{\rm id}_Y$,
defined on the set $\widehat{\mathcal K^\tau}$ of all
biholomorphic self-mappings of $\mathbb{D}_1$ with the boundary
attractive fixed point $\tau\in\pl\mathbb{D}_1$. As above,
condition (i) follows by the selection of an appropriate branch of
the power, conditions (ii) and (iii) hold automatically.
Furthermore, by a multidimensional analog of the boundary
Wolff--Schwarz Lemma (see, for example, \cite{RS-SD-book})
\[
\frac{1-\|x\|_X^2}{|1-\langle x,\tau\rangle|^2}\le
\frac{1-\|f(x)\|_X^2}{|1-\langle f(x), \tau\rangle|^2}\,.
\]
Therefore,
\[
\left\|\widehat{\Gamma}(f,x)\right\|_{L(Y)}= \left|\frac{1-\langle
f(x), \tau\rangle}{1-\langle x,\tau\rangle}\right|^{2/r} \le
\frac{(1-\|f(x)\|_X^2)^{1/r}}{(1-\|x\|_X^2)^{1/r}}\,,
\]
i.e., condition (iv) is satisfied.
\eexa

For each appropriate mapping $\widehat{\Gamma}$, one corresponds
the extension operator $\widehat\Phi:\widehat{\mathcal
K}\mapsto\Hol(\mathbb{D})$ defined by
\begin{equation}\label{RS-abs1}
\widehat\Phi[f](x,y)=\widehat\Phi_{\widehat{\Gamma}}[f](x,y)=\left(f(x),\widehat{\Gamma}(f,x)y
\right).
\end{equation}

In Section~\ref{sect-semi} below, we will study its modification
as an extension operator for one-parameter semigroups.

\begin{lemma}\label{lem-1}
Let $\widehat{\Gamma}:\,\widehat{\mathcal
K}\times\mathbb{D}_1\mapsto L(Y)$ be appropriate. Let $f,g\in
\widehat{\mathcal K}$. Then

(a) $\widehat\Phi[f]\in\Hol(\mathbb{D})$;

(b) $\widehat\Phi[f\circ g]=\widehat\Phi[f]\circ\widehat\Phi[g]$.
\end{lemma}

For the original Roper--Suffridge operator (\ref{RS}), assertion
(a) of this lemma can be found in \cite{B-K}.

\pr Assertion (a) means that for each point $(x,y)\in\mathbb{D}$
the inequality
\[
\left\|\widehat{\Gamma}(f,x)y \right\|_Y <
p\left(\left\|f(x)\right\|_X \right)
\]
holds. Indeed, since $(x,y)\in\mathbb{D}$, we have
$\|y\|_Y<p(\|x\|_X)$. Therefore,
\[
\left\|\widehat{\Gamma}(f,x)y \right\|_Y \le \left\|
\widehat{\Gamma}(f,x) \right\|_{L(Y)} \|y\|_Y \le \frac{p
\left(\|f(x)\|_X\right)}{p\left(\|x\|_X\right)}\,\|y\|_Y < p
\left(\|f(x)\|_X\right).
\]

To prove assertion (b), we just calculate:
\begin{eqnarray*}
&&\biggl(\widehat\Phi[f]\circ\widehat\Phi[g]\biggr)(x,y) =
\widehat\Phi[f]\biggl(
\widehat\Phi[g](x,y) \biggr) =\widehat\Phi[f]\biggl(g(x),\widehat{\Gamma}(g,x)y\biggr) \\
&&= \biggl(f(g(x)),\widehat{\Gamma}(f,g(x))
\widehat{\Gamma}(g,x)y\biggr)
=\biggl((f\circ g)(x),\widehat{\Gamma}(f\circ g,x)y\biggr) \\
&&=\biggl(\widehat\Phi[f\circ g]\biggr)(x,y).
\end{eqnarray*}
\epr

Now we expand the notion of appropriate operators to biholomorphic
mappings $\mathbb{D}_1\mapsto X.$

\begin{defin}\label{def-approp-a}
Let a set $\widehat{\mathcal K}\subset\Hol(\mathbb{D}_1)$ and an
appropriate mapping $\widehat{\Gamma}$ be given. Suppose that
there are (a) a non-empty set
$\mathcal{K}=\mathcal{K}_{\mathbb{D}_1}\subset\Hol(\mathbb{D}_1,X)$
consisting of biholomorphic mappings and (b) a mapping
$\Gamma=\Gamma_{\mathbb{D}_1}:\mathcal{K}\times\mathbb{D}_1
\mapsto L(Y)$ continuous on $\mathcal{K}$ and holomorphic on
$\mathbb{D}_1$ such that

\begin{itemize}
\item[(i)] for all $h_1,h_2\in\mathcal{K}$ with
$h_1(\mathbb{D}_1)\subset h_2(\mathbb{D}_1)$, we have
$h_2^{-1}\circ h_1\in\widehat{\mathcal K}$;

\item[(ii)] $\Gamma(h,g(x))\widehat{\Gamma}(g,x)= \Gamma(h\circ
g,x)$ for all $h\in\mathcal{K}$, $g\in\widehat{\mathcal K}$ and
$x\in\mathbb{D}_1$;

\item[(iii)] for each $h\in\mathcal{K}$ and $x\in\mathbb{D}_1$,
the operator $\Gamma(h,x)$ is invertible.
\end{itemize}
Then we say that $\Gamma=\Gamma_{\mathbb{D}_1}$ is appropriate.
\end{defin}

\rema\label{rem8} For appropriate mappings considered in
Examples~\ref{ex1} and~\ref{ex2} and defined on the set of
mappings normalized by $f(0)=0$, one can choose $\Gamma$ to be
defined by the same formula as $\widehat{\Gamma}$, that is,
respectively, $\Gamma(h,x)=\left(J_h(x)\right)^\alpha{\rm id}_Y$
or $\Gamma(h,x)=\left(\frac{h(x)}{x}\right)^\beta{\rm id}_Y$,
where $h(0)=0$. As previously mentioned, the operator
$\widehat{\Gamma}$ from Example~\ref{ex1} can be defined on the
set $\widehat{\mathcal{K}^\tau},\ \tau\in\overline{\mathbb{D}_1}$.
In this case, we can again use the same formula.

Concerning Example~\ref{ex3}, it is possible to proceed as
follows. We choose some mapping $A:X\mapsto L(Y)$ and a set of
biholomorphic mappings $\mathcal{K}$ such that $A(h(x))$ is
invertible for all $h\in\mathcal{K}$ and $x\in\mathbb{D}_1$. Then
we set $\displaystyle\Gamma(h,x)=(1-\langle x,\tau\rangle)
^{-2/r}A(h(x))$. For instance, we can choose $\mathcal{K}$ to be a
set of biholomorphic mappings $h\in\Hol(\mathbb{D}_1,X)$ with
$\langle h(x),\tau\rangle\not=0$ for all $x\in\mathbb{D}_1$ and to
define $\Gamma(h,x)=\displaystyle\left(\frac{\langle
h(x),\tau\rangle} {1-\langle x, \tau \rangle}\right)^{2/r}{\rm
id}_Y$. \erema

Similar to (\ref{RS-abs1}), we define the extension operator $\Phi
:\mathcal{K}\mapsto\Hol(\mathbb{D},Z)$~by
\begin{equation}\label{RS-abs}
\Phi[h](x,y)=\Phi_\Gamma[h](x,y)=\left(h(x),\Gamma(h,x)y \right).
\end{equation}

This operator will be the main subject in Section~\ref{sect-star}.
In particular, we will study its action on starlike and spirallike
mappings.

\begin{lemma}\label{lem-2}
Let $\widehat{\Gamma}:\,\widehat{\mathcal
K}\times\mathbb{D}_1\mapsto L(Y)$ and
$\Gamma:\,\mathcal{K}\times\mathbb{D}_1\mapsto L(Y)$ be
appropriate. Let $h\in\mathcal{K}$ and $g\in\widehat{\mathcal K}$.
Then
\[
\Phi[h\circ g]=\Phi[h]\circ \widehat\Phi[g].
\]
In addition, $\Phi[h]$ is biholomorphic, and for
$(z,w)\in\Phi[h](\mathbb{D})$ we have
\[
\left(\Phi[h]\right)^{-1}(z,w)=\left(h^{-1}(z),\left(
\Gamma(h,h^{-1}(z)) \right)^{-1}w\right).
\]
\end{lemma}
\pr The first assertion follows by the calculation:
\begin{eqnarray*}
&&\biggl(\Phi[h]\circ\widehat\Phi[g]\biggr)(x,y) =
\Phi[h] \biggl(\widehat\Phi[g](x,y) \biggr) = \Phi[h]\biggl(g(x),\widehat{\Gamma}(g,x)y\biggr) \\
&& = \biggl(h(g(x)),\Gamma(h,g(x))\widehat{\Gamma}(g,x)y\biggr)=
\biggl((h\circ g)(x),\Gamma(h\circ g,x)y\biggr) \\
&&=\biggl(\Phi[h\circ g]\biggr)(x,y).
\end{eqnarray*}

The last assertion is obvious. \epr

Now we are ready to turn to appropriate operators on domains
biholomorphically equivalent to the unit ball $\mathbb{D}_1$.

\begin{defin}\label{def-approp1}
Let $\Gamma:\mathcal{K}\times\mathbb{D}_1\mapsto L(Y)$ be
appropriate. Given a domain $\Omega\in X$ biholomorphically
equivalent to the ball $\mathbb{D}_1$, we define the set
$\mathcal{K}_\Omega$ to consist of all biholomorphic mappings
$f\in\Hol(\Omega,X)$ for which there is a biholomorphic mapping
$h$ of $\,\mathbb{D}_1$ onto $\Omega$ such that both $h$ and
$f\circ h$ belong to $\mathcal{K}$. For $f\in\mathcal{K}_\Omega$
and $x\in\Omega$ we define the appropriate mapping $\Gamma_\Omega$
by
\[
\Gamma_\Omega(f,x):=\Gamma(f\circ h,h^{-1}(x))\left(\Gamma(h,
h^{-1}(x))\right)^{-1}.
\]
\end{defin}

The next assertion can be checked directly.
\begin{lemma}\label{lem-3}
The mapping $\Gamma_\Omega$ is well-defined in the sense that it
is independent of the choice of a biholomorphic mapping
$h\in\mathcal{K}$ of $\mathbb{D}_1$ onto $\Omega$. Moreover,
$\Gamma_\Omega$ has the following properties:
\begin{itemize}
\item[(i)] $\Gamma_\Omega({\rm id}_X,x)={\rm id}_Y$ for all
$x\in\Omega$;

\item[(ii)]
$\Gamma_\Omega(f,g(x))\Gamma_\Omega(g,x)=\Gamma_\Omega(f\circ
g,x)$ for all $f\in\mathcal{K}_\Omega$,
$g\in\mathcal{K}_\Omega\cap\Hol(\Omega)$ and $x\in\Omega$;

\item[(iii)] for each $f\in\mathcal{K}_\Omega$ and $x\in\Omega$,
the operator $\Gamma_\Omega(f,x)$ is invertible. In particular, if
$h\in\mathcal{K}$ is a biholomorphic mapping of $\mathbb{D}_1$
onto $\Omega$, then $h^{-1}\in\mathcal{K}_\Omega$ and
$\Gamma_\Omega(h^{-1},h(x))= \left(\Gamma(h,x) \right)^{-1}$
\end{itemize}
\end{lemma}

\pr Let $h_1,h_2$ be biholomorphic mappings of $\mathbb{D}_1$ onto
$\Omega$ such that\newline ${h_1,h_2,f\circ h_1,f\circ
h_2}\in\mathcal{K}$. We have to show that
\[
\Gamma(f\circ h_1,h_1^{-1}(x))\left(\Gamma(h_1,
h_1^{-1}(x))\right)^{-1} = \Gamma(f\circ
h_2,h_2^{-1}(x))\left(\Gamma(h_2, h_2^{-1}(x))\right)^{-1}.
\]
for all $x\in\Omega$. Denote $w=h_1^{-1}(x)$ and let
$\phi:=h_2^{-1}\circ h_1$ be an automorphism of $\mathbb{D}_1$
which belongs to $\mathcal{K}$ by Definition~\ref{def-approp-a}.
Then the equality above can be rewritten as
\[
\Gamma(f\circ h_1,w)\left(\Gamma(h_1,w)\right)^{-1} =
\Gamma(f\circ h_2,\phi(w))\left(\Gamma(h_2, \phi(w))\right)^{-1}.
\]
This relation holds by Definition~\ref{def-approp-a} since
$h_2\circ\phi=h_1$.

Properties (i) and (iii) hold by Definition~\ref{def-approp1}. To
check property~(ii), let consider the expression
\begin{eqnarray*}
&&\Gamma_\Omega(f,g(x))\Gamma_\Omega(g,x)\Gamma(h,h^{-1}(x)) \\&=&
\Gamma(f\circ h,h^{-1}(g(x)))\left(\Gamma(h,h^{-1}(g(x)))
\right)^{-1}\Gamma(g\circ h,h^{-1}(x)) \\&=& \Gamma(f\circ
h,\psi(w))\left(\Gamma(h,\psi(w) \right)^{-1}\Gamma(h\circ\psi,w),
\end{eqnarray*}
where we denote $\psi=h^{-1}\circ g\circ h$ and $w=h^{-1}(x)$.
Since $(g\circ h)(\mathbb{D}_1)\subset h(\mathbb{D}_1)$, we
conclude by Definition~\ref{def-approp-a}~(i) that
$\psi\in\widehat{\mathcal{K}}$. Now, using condition (ii) of
Definition~\ref{def-approp-a}, we obtain:
\begin{eqnarray*}
&&\Gamma_\Omega(f,g(x))\Gamma_\Omega(g,x)\Gamma(h,h^{-1}(x)) \\&=&
\Gamma(f\circ h,\psi(w))\widehat{\Gamma}(\psi,w) = \Gamma(f\circ
h\circ\psi,w) = \Gamma(f\circ g\circ h,h^{-1}(x)),
\end{eqnarray*}
so (ii) follows.
\epr

{\bf In what follows, all operator-valued mappings
$\widehat{\Gamma},\ \Gamma$ and $\Gamma_\Omega$ are assumed to be
appropriate.}

\section{Extension operators for semigroups}\label{sect-semi}
\setcounter{equation}{0}

In this section we study extension operators for one-parameter
continuous semigroups. It turns out that for a given appropriate
mapping, each semigroup on the unit ball of $X$ admits a family of
extensions.

\begin{theorem}\label{th-sg}
Let $\widehat{\Gamma}:\,\widehat{\mathcal
K}\times\mathbb{D}_1\mapsto L(Y)$ be appropriate, i.e., conditions
(i)--(iv) of Definition~\ref{def-approp} are satisfied. Let
$S=\{F_t\}_{t\ge0}\subset\Hol(\mathbb{D}_1)$ be a semigroup on the
ball $\mathbb{D}_1$ such that $S\subset\widehat{\mathcal K}$. Let
$\Sigma=\{G_t\}_{t\ge0}$ be a semigroup on the ball $\mathbb{D}_2$
such that each its element $G_s,\ s\ge0,$ satisfies
$\|G_s(y)\|_Y\le\|y\|_Y$ for all $y\in\mathbb{D}_2$ and commutes
with operators $\widehat{\Gamma}(F_t,x)$ for all $t\ge0$ and
$x\in\mathbb{D}_1$:
\begin{equation}\label{comm}
\widehat{\Gamma}(F_t,x)\circ G_s= G_s\circ\widehat{\Gamma}(F_t,x).
\end{equation}
Then the family $\widetilde{S}=\left\{\widetilde{F_t}
\right\}_{t\ge0}$ defined by
\begin{equation}\label{RS-sg-1}
\widetilde{F_t}(x,y)=\biggl(F_t(x),\widehat{\Gamma}(F_t,x)\,
G_t(y)\biggr),
\end{equation}
forms a semigroup on $\mathbb{D}$.
\end{theorem}
\rema\label{rem4} In the case when $\Sigma$ is a uniformly
continuous semigroup of proper contractions (hence, $G_t=e^{-Bt}$
for some accretive operator $B$, see \cite{RS-SD-book}), the
commutativity condition (\ref{comm}) can be replaced by the
following one: {\it all operators $\widehat{\Gamma}(F_t,x),\
t\ge0,\ x\in\mathbb{D}_1,$ commute with $B$.} In particular, the
last condition always holds if $\ \widehat{\Gamma}(F_t,x)$ is a
scalar operator for each $t\ge0$ and $x\in\mathbb{D}_1$. \erema

\pr Since $G_t$ is a contraction, it follows by
Lemma~\ref{lem-1}~(a) that $\widetilde{F_t}$ is a self-mapping of
$\mathbb{D}$ for each $t\ge0$. The continuity of
$\widehat{\Gamma}$ and condition (i) of
Definition~\ref{def-approp} imply that
\begin{eqnarray*}
\lim_{t\to 0^+}\widetilde{F_t}(x,y) =\lim_{t\to 0^+}\left(
F_t(x),\widehat{\Gamma}(F_t,x)\, G_t(y)\right) =
\left(x,\widehat{\Gamma}({\rm id_X}, x)\,G_0( y))\right) =(x,y).
\end{eqnarray*}
Similarly to the proof of Lemma~\ref{lem-1}~(b), we have for all
$t,s>0$:
\begin{eqnarray*}
&& \biggl(\widetilde{F_t}\circ\widetilde{F_s}\biggr)(x,y)
=\widetilde{F_t}\biggl(\widetilde{F_s}(x,y)\biggr)
=\widetilde{F_t}\biggl(F_s(x),\widehat{\Gamma}(F_s,x)\, G_s(y)\biggr) \\
&&= \biggl(F_t(F_s(x)),\widehat{\Gamma}(F_t,F_s(x))\circ
G_t\circ\widehat{\Gamma}(F_s,x)\circ G_s(y)\biggr)
\\
&&=
\biggl(F_t(F_s(x)),\widehat{\Gamma}(F_t,F_s(x))\,\widehat{\Gamma}(F_s,x)\circ
G_t\circ G_s(y)\biggr)
\\&&=
\biggl((F_t\circ F_s)(x),\widehat{\Gamma}(F_t\circ
F_s,x)\,G_{t+s}(y)\biggr)=\widetilde{F_{t+s}}(x,y).
\end{eqnarray*}
This calculation completes the proof. \epr

\begin{corol}\label{cor1-fp}
Let an appropriate mapping $\widehat{\Gamma}$ and semigroups
$S\subset\widehat{\mathcal{K}}$ and $\Sigma\subset
\Hol(\mathbb{D}_2)$ be as above. Denote by
$\mathcal{M}\subset\mathbb{D}_1$ the stationary point set of $S$.
Then the stationary point set $\widetilde{\mathcal{M}}$ of the
extended semigroup $\widetilde S$ satisfies the following
inclusion:
\[
\left\{(x,0)\in\mathbb{D}: x\in\mathcal{M} \right\} \subset
\widetilde{\mathcal{M}} \subset \left\{(x,y)\in\mathbb{D}:
x\in\mathcal{M} \right\}.
\]
\end{corol}

To find the semigroup generator, we require the Frech\'et
differentiability of $\widehat{\Gamma}$ in $f\in\widehat{\mathcal
K}$, namely, {\it
\begin{itemize}
\item at each point $f\in\widehat{\mathcal K}$ the Frech\'et
derivative (denoted by $\partial\widehat{\Gamma}(f,x)$) exists as
a linear operator defined on ${\rm span}(\widehat{\mathcal K})$.
\end{itemize}}

Just differentiating (\ref{RS-sg-1}) at $t=0^+$, we obtain the
following assertion.

\begin{corol}\label{cor1-gen}
Let an appropriate mapping $\widehat{\Gamma}$ and semigroups
$S\subset\widehat{\mathcal K}$ and $\Sigma\subset
\Hol(\mathbb{D}_2)$ be as above, and let condition ($\bullet$) be
satisfied. If $S$ is generated by a mapping
$f\in\Hol(\mathbb{D}_1,X)$, and $\Sigma$ is generated by a mapping
$g\in\Hol(\mathbb{D}_2,Y)$, then the extended semigroup
$\widetilde S$ defined by (\ref{RS-sg-1}) is generated as well.
Its generator $\widetilde f$ is defined by
\[
\widetilde{f}(x,y)=\biggl(f(x),\partial\widehat{\Gamma}({\rm
id}_X,x)[f]y+ g(y)\biggr).
\]
\end{corol}

We proceed with the extension of conjugate semigroups.

\begin{theorem}\label{th-conj}
Let $\{F_t\}_{t\ge0}\subset\widehat{\mathcal K}$ and
$\{\Psi_t\}_{t\ge0}\subset\mathcal{K}_\Omega$ be conjugate
semigroups acting on the unit ball $\mathbb{D}_1$ and a domain
$\Omega\subset X$, respectively. Let
$h\in\Hol(\mathbb{D}_1,\Omega)\cap\mathcal{K}$ be their
intertwining map. Then the mapping $\widetilde{h}=\Phi[h]$ defined
by~(\ref{RS-abs}) is the intertwining map for the semigroup
$\widetilde{S}=\left\{\widetilde{F_t}\right\}_{t\ge0}$ defined
by~(\ref{RS-sg-1}) and the semigroup $\left\{\widetilde{\Psi_t}
\right\}_{t\ge0}$ acting on $\Phi[h](\mathbb{D})$ and defined by
\begin{equation} \label{eq-psi}
\widetilde{\Psi_t}(z,w)=\left(\Psi_t(z),\Gamma_\Omega(\Psi_t,z)\tilde{G}_t(z,w)
\right),
\end{equation}
where
\[
\widetilde{G}_t(z,w)=  \Gamma(h, h^{-1}(z)) G_t
\left(\Gamma_\Omega(h^{-1},z)w\right).
\]
\end{theorem}

Note that if all mappings $G_t, t\ge 0,$ commute with
$\Gamma(h,x)$ (for example, in the case described in
Remark~\ref{rem4}), then $\widetilde{G}_t(z,w)=G_t(w)$.

\pr It has already be proven in Theorem~\ref{th-sg} that the
family $\widetilde{S}=\left\{\widetilde{F_t}\right\}_{t\ge0}$
forms a semigroup on $\mathbb{D}$. Therefore, the family
$\left\{\Phi[h]\circ\widetilde{F_t}\circ\left(\Phi[h]\right)^{-1}
\right\}_{t\ge0}$ forms a semigroup on $\Phi[h](\mathbb{D})$ which
is conjugate to $\widetilde S$ with the intertwining mapping
$\Phi[h]$. Let us find its exact form. By Lemmas~\ref{lem-2}
and~\ref{lem-3},
\[
\left(\Phi[h]\right)^{-1}(z,w)=\left(h^{-1}(z),\left(
\Gamma(h,h^{-1}(z)) \right)^{-1}w\right)=
\left(h^{-1}(z),\Gamma_\Omega(h^{-1},z)w\right).
\]
Now, we substitute
\begin{eqnarray*}
&& \widetilde{F_t}\circ\left(\Phi[h]\right)^{-1}(z,w) =
\widetilde{F_t} \left(h^{-1}(z),\Gamma_\Omega(h^{-1},z)w\right)
\\ && = \left(F_t\left(h^{-1}(z)\right),
\widehat\Gamma\left(F_t,h^{-1}(z)\right) G_t
\left(\Gamma_\Omega(h^{-1},z)w\right) \right).
\end{eqnarray*}
By Definition~\ref{def-approp-a},
$\widehat\Gamma\left(F_t,h^{-1}(z)\right)=\left(\Gamma(h,F_t\circ
h^{-1}(z)) \right)^{-1}\Gamma(h\circ F_t,h^{-1}(z))$. In addition,
since $h$ is the intertwining map for $\{F_t\}_{t\ge0}$ and
$\{\Psi_t\}_{t\ge0}$, we conclude that $F_t\circ
h^{-1}=h^{-1}\circ\Psi_t$. Therefore,
\begin{eqnarray*}
&&\hspace{-5mm}  \left(\Phi[h]\circ\widetilde{F_t}\circ
\left(\Phi[h]\right)^{-1}\right)(z,w)
\\ &&\hspace{-5mm}  = \Phi[h]\left(
h^{-1}\left(\Psi_t(z)\right),\left(\Gamma(h,h^{-1}\circ \Psi_t(z))
\right)^{-1}\Gamma(\Psi_t\circ h,h^{-1}(z)) G_t
\left(\Gamma_\Omega(h^{-1},z)w\right)                \right) \\
&&\hspace{-5mm} = \left( \Psi_t(z),\Gamma(\Psi_t\circ h,h^{-1}(z))
G_t \left(\Gamma_\Omega(h^{-1},z)w\right)                \right).
\end{eqnarray*}
Finally, by Definition~\ref{def-approp1}
\[
\Gamma(\Psi_t\circ h,h^{-1}(z))=\Gamma_\Omega(\Psi_t,z)\Gamma(h,
h^{-1}(z)).
\]
Thus,
\begin{eqnarray*}
&& \left(\Phi[h]\circ\widetilde{F_t}\circ
\left(\Phi[h]\right)^{-1}\right)(z,w) \\&&= \left(
\Psi_t(z),\Gamma_\Omega(\Psi_t,z)\Gamma(h, h^{-1}(z)) G_t
\left(\Gamma_\Omega(h^{-1},z)w\right)                \right),
\end{eqnarray*}
and the assertion follows. \epr

\section{Starlikeness, spirallikeness and convexity in one
direction}\label{sect-star}

\setcounter{equation}{0}

The main results of this section are Theorems~\ref{th-spiral}
and~\ref{th-conv_dir}  below. In these theorems, given a
biholomorphic mapping $h$, we examine geometric properties of its
extension $\Phi[h]$ defined by formula~(\ref{RS-abs}):
$\Phi[h](x,y)= \left(h(x),\Gamma(h,x)y \right)$.

\begin{theorem}\label{th-spiral}
Let $h\in\Hol(\mathbb{D}_1,X)$ be an $A$-spirallike mapping.
Suppose that $e^{-At}\circ h\in\mathcal{K}$ for all $t\ge0$ and
there is $C\in L(Y)$ such that
\begin{equation}\label{assump}
\Gamma\left(e^{-At}\circ h,x\right) =
e^{-Ct}\Gamma(h,x)\quad\mbox{for all}\quad t\ge0\mbox{ and
}x\in\mathbb{D}_1.
\end{equation}
Then the mapping $\Phi[h]$ is $\left(
\begin{array}{cc}
  A & 0 \\
  0 & B+C \\
\end{array}%
\right) $-spirallike for any accretive operator $B\in L(Y)$ which
commutes with $C$ and with $\Gamma(h,x)$ for all
$x\in\mathbb{D}_1$ and such that the function $\Re\lambda$ is
bounded away from zero on the spectrum of $B+C$.
\end{theorem}

By definition, for each point of the image of a spirallike mapping
there is a spiral curve which is contained in the image. Our
theorem asserts that the image of the extension of a spirallike
mapping contains not only a one-dimensional spiral curve but at
least some manifold (of real codimension $(2n-1)$ when $X=\C^n$).
We illustrate this effect in Examples~\ref{ex4} and~\ref{ex5}
below.

\begin{theorem}\label{th-conv_dir}
Let a biholomorphic mapping $h\in\Hol(\mathbb{D}_1,X)$ be convex
in the direction $\tau$, where $\tau\in\pl\mathbb{D}_1$. Suppose
that $h+t\tau\in\mathcal{K}$ for all $t\ge0$ and there is $C\in
L(Y)$ such that
\begin{equation}\label{assump-c}
\Gamma\left(h+t\tau,x\right) = e^{-Ct}\Gamma(h,x)\quad\mbox{for
all}\quad t\ge0\mbox{ and }x\in\mathbb{D}_1.
\end{equation}
Then for each point $(z,w)\in\Phi[h](\mathbb{D})$, the set
$\Phi[h](\mathbb{D})$ contains the curve
\[
\biggl\{\left(z+t\tau,e^{-(B+C)t}w \right),\ t\ge0 \biggr\},
\]
for any accretive operator $B\in L(Y)$ which commutes with $C$ and
with $\Gamma(h,x)$ for all $x\in\mathbb{D}_1$ and such that the
function $\Re\lambda$ is non-negative on the spectrum of $B+C$. In
particular, if $\Gamma\left(h+t\tau,x\right) =\Gamma(h,x)$ for all
$t\ge0$ and $x\in\mathbb{D}_1$, then the mapping $\Phi[h]$ is
convex in the direction $(\tau,0)$.
\end{theorem}

\exa\label{ex4} Let $X=\C^n$ with an arbitrary norm. Similar to
the examples in Section~\ref{sect-approp}, we define the unit ball
in the space $Z=X\times Y=\mathbb{C}\times Y$ by
\[
\mathbb{D}=\left\{(x,y):\ \|x\|_X^2+ \|y\|_Y^r<1 \right\},\qquad
r\ge1.
\]
Consider the appropriate mapping
$\Gamma(h,x)=\left(J_h(x)\right)^{\frac2{r(n+1)}}{\rm id}_Y$ (cf.,
Example~\ref{ex1} and Remark~\ref{rem8} above) and the
corresponding extension operator
\[
\Phi[h](x,y)=\left( h(x), \left(J_h(x)
\right)^{\frac2{r(n+1)}}y\right).
\]

(1) Let $A$ be a diagonal matrix, $A={\rm
diag}(\mu_1,\ldots,\mu_n)$ with $\Re\mu_j>0$. Take any
$A$-spirallike mapping $h$ on the unit ball $\mathbb{D}_1\subset
X$ with respect to either an interior or a boundary point. Since
$J_{e^{-At}h}(x)=e^{-\tr A \cdot t}J_h(x)$, where $\tr
A=\mu_1+\ldots+\mu_n$ is the trace of the matrix $A$, we get that
the operator $C$ in formula~(\ref{assump}) is given by
$C=\frac{2\tr A}{r(n+1)}\,{\rm id}_Y$. According to
Theorem~\ref{th-spiral}, the extended mapping $\Phi[h]$ is $\left(
\begin{array}{cc}
  A & 0 \\
  0 & B+\frac{2\tr A}{r(n+1)}\,{\rm id}_Y \\
\end{array}\right)$-spirallike for any accretive operator $B\in
L(Y)$. To understand this effect, consider the simplest case
$Y=\C$. In this situation all linear operators are just
multiplication by scalars. We have that for any point
$(z_0,w_0)\in\Phi[h](\mathbb{D})$, the image $\Phi[h](\mathbb{D})$
contains the set
\[
\left\{(z,w):\ z=e^{-A t}z_0,\ w=e^{-\left(\lambda+\frac{\tr A}
r\right)t} w_0,\ t\ge0,\ \Re\lambda>0 \right\},
\]
or, equivalently,
\[
\left\{(z,w):\ z=e^{-A t}z_0,\ |w|<e^{-\frac{t\Re(\tr A)}r}
\left|w_0\right|,\ t\ge0 \right\}.
\]
Schematically, this set is presented in Fig.~1.

\begin{center}
\begin{figure}[h]
\centering
\includegraphics[angle=0,width=8.5cm,totalheight=5.5cm]{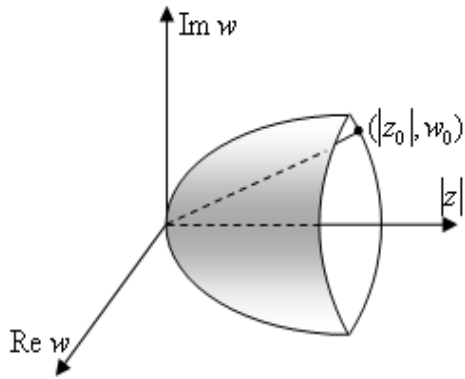}
\caption{The `spiral' segment and the manifold inside
$\Phi[h](\mathbb{D})$.}
\end{figure}
\end{center}

(2) Let now $\tau\in\partial\mathbb{D}_1$. Take a mapping $h$
convex in the direction $\tau$. Since $J_{h+t\tau}(x)=J_h(x)$, we
conclude that the operator $C$ in formula~(\ref{assump-c}) is
zero. According to Theorem~\ref{th-conv_dir}, for each point
$(z_0,w_0)\in\Phi[h](\mathbb{D})$, the image $\Phi[h](\mathbb{D})$
contains the set
\[
\biggl\{\left(z_0+t\tau,e^{-Bt}w_0 \right),\ t\ge0, \ B\in L(Y)\
\mbox{is attractive} \biggr\}.
\]
Once again, we restrict our consideration to the case $Y=\C$. Then
for any point $(z_0,w_0)\in\Phi[h](\mathbb{D})$, the image
$\Phi[h](\mathbb{D})$ contains the set
\[
\left\{(z,w):\ z=z_0+t\tau,\ |w|\le \left|w_0\right|,\ t\ge0
\right\}
\]
(see Fig. 2). \eexa

\begin{center}
\begin{figure}[h]
\centering
\includegraphics[angle=0,width=8.5cm,totalheight=5cm]{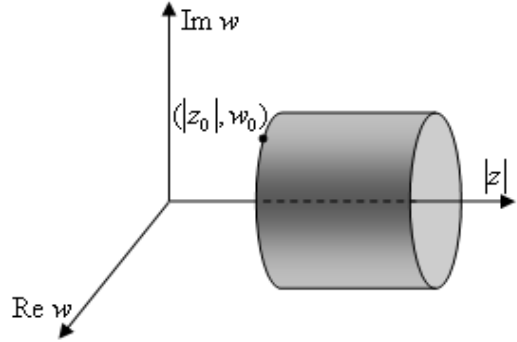}
\caption{The cylindric set inside $\Phi[h](\mathbb{D})$.}
\end{figure}
\end{center}

\exa\label{ex5} Let $X$ be a complex Hilbert space. As above, we
define the unit ball in the space $Z=X\times Y$ by
\[
\mathbb{D}=\left\{(x,y):\ \|x\|_X^2+ \|y\|_Y^r<1 \right\},\qquad
r\ge1.
\]

Let $A\in L(X)$ be a bounded linear operator such that the
function $\Re\lambda$ is bounded away from zero on the spectrum of
$A$. Suppose that a vector $\tau\in X,\ \|\tau\|_X=1,$ is an
eigenvector of the adjoint operator $A^*$, and $\bar\lambda$ is
the corresponding eigenvalue.

Consider the appropriate mapping
$\Gamma(h,x)=\displaystyle\left(\frac{\langle h(x),\tau\rangle}
{1-\langle x, \tau \rangle}\right)^{2/r}{\rm id}_Y$ (cf.,
Example~\ref{ex3} and Remark~\ref{rem8} above) and the
corresponding extension operator
\[
\Phi[h](x,y)=\left( h(x), \left(\frac{\langle h(x),\tau\rangle}
{1-\langle x, \tau \rangle}\right)^{2/r}y\right).
\]

Take an $A$-spirallike mapping $h\in\Hol(\mathbb{D}_1,X)$ with
respect to a boundary point with $\lim\limits_{r\to1^-}
h(r\tau)=0$ and such that $\langle h(x),\tau\rangle\not=0$ for all
$x\in\mathbb{D}_1$. Clearly,
\[
\langle e^{-At}h(x),\tau\rangle=\langle h(x),e^{-A^*t}\tau\rangle
=e^{-\lambda t}\langle h(x),\tau\rangle.
\]
Therefore, equality~(\ref{assump}) holds for the operator
$C=\frac{2\lambda}{r}\,{\rm id}_Y$. So, Theorem~\ref{th-spiral}
asserts that the extended mapping $\Phi[h]$ is $\left(
\begin{array}{cc}
  A & 0 \\
  0 & B+\frac{2\lambda}{r}\,{\rm id}_Y \\
\end{array}\right)$-spirallike for any accretive operator $B\in
L(Y)$.

In the particular case when $X$ and $Y$ are one-dimensional, we
conclude that for any $\lambda,\ \Re\lambda>0,$ the extension of
each $\lambda$-spirallike function with respect to a boundary
point is $\left(
\begin{array}{cc}
  \lambda & 0 \\
  0 & \mu+\frac{2\lambda}{r} \\
\end{array}\right)$-spirallike for any number $\mu$ with non-negative real
part. We see that if $r<2$ then the extended mapping may be not
$\lambda$-spirallike. \eexa

As previously mentioned in Section~\ref{sect-nota}, the images of
spirallike mappings and mappings convex in one direction are
invariant under action of a linear semigroup of proper
contractions and a semigroup of shifts, respectively. More
generally, we can consider a semigroup of affine mappings. Thus,
both Theorems~\ref{th-spiral} and \ref{th-conv_dir} can be
considered as consequences of the following general assertion,
where we denote by $\Sigma=\Sigma(A,\lambda,\tau)=\left\{\Psi_t
\right\}_{t\ge0}$ the semigroup of affine mappings defined by
\[
\Psi_t(z)=e^{-At}z +\lambda\int_0^te^{-As}\tau ds,
\]
where $A\in L(X),\ \lambda\ge0$ and $\tau\in X,\ \|\tau\|_X=1$.

\begin{theorem}\label{th-spiral-gen}
Let $\Sigma=\Sigma(A,\lambda,\tau)$ be a semigroup of affine
mappings. Let $h\in\Hol(\mathbb{D}_1,X)$ be biholomorphic, and
$h(\mathbb{D}_1)$ be $\Sigma$-invariant. Suppose that $\Psi_t\circ
h\in\mathcal{K}$ for all $t\ge0$ and there is an operator $C\in
L(Y)$ such that
\begin{equation}\label{assump-1}
\Gamma\left(\Psi_t\circ h,x\right) = e^{-Ct}\Gamma(h,x)
\end{equation}
for all $t\ge0$ and $x\in\mathbb{D}_1$. Let
$\left\{G_s\right\}_{s\ge0}\subset\Hol(\mathbb{D}_2)$ be a
semigroup such that each its element $G_s,\ s\ge0,$ satisfies
$\|G_s(y)\|_Y\le\|y\|_Y$ for all $y\in\mathbb{D}_2$ and commutes
with $\Gamma(h,x)$ for all $x\in\mathbb{D}_1$.

Then for each point $(z,w)\in\Phi[h](\mathbb{D})$, the image
$\Phi[h](\mathbb{D})$ contains the curve
\[
\biggl\{\left(\Psi_t(z),e^{-Ct}G_t(w) \right),\ t\ge0 \biggr\}.
\]
\end{theorem}

\pr Since $h(\mathbb{D}_1)$ is $\Sigma$-invariant, the family
$S=\left\{F_t\right\}_{t\ge0}$ with $F_t(x)=h^{-1}\circ
\Psi_t\circ h(x)$ forms a semigroup on $\mathbb{D}_1$. By our
assumption and condition~(i) of Definition~\ref{def-approp-a} we
conclude that $S\subset\widehat{\mathcal{K}}$. Obviously, $h$ is
the intertwining map for the semigroups $S$ and $\Sigma$ (acting
on the domain $\Omega=h(\mathbb{D}_1)$).

By Theorem~\ref{th-conj}, the image of the mapping $\Phi[h]$
contains together with each point $(z,w)\in\Phi[h](\mathbb{D})$
the whole semigroup trajectory $\left\{\widetilde{\Psi_t}(z,w),\
t\ge0\right\}$, where $\widetilde{\Psi_t}$ is defined by
(\ref{eq-psi}). It follows by (\ref{assump-1}) that
\begin{eqnarray*}
\Gamma_\Omega(\Psi_t,z)= \Gamma(\Psi_t\circ
h,h^{-1}(z))\left(\Gamma(h,
h^{-1}(z))\right)^{-1}  \\
=e^{-Ct}\Gamma(h,h^{-1}(z))\left(\Gamma(h,
h^{-1}(z))\right)^{-1}=e^{-Ct}.
\end{eqnarray*}
In addition, $\widetilde{G}_t(z,w)=G_t(w)$. Therefore,
\[
\widetilde{\Psi_t}(z,w)=\left(\Psi_t(z),e^{-Ct}G_t(w) \right),
\]
and the assertion is proved. \epr

\noindent{\bf Proof of Theorem~\ref{th-spiral}.} Let
$h\in\Hol(\mathbb{D}_1,X)$ be an $A$-spirallike mapping which
satisfies~(\ref{assump}). Let $B\in L(Y)$ be an accretive operator
which commutes with $C$. Then the semigroup
$\left\{e^{-Bs}\right\}_{s\ge0}$ consists of proper contractions
with respect to the norm $\|\cdot\|_Y$. Since
\begin{eqnarray*}
\widehat{\Gamma}(F_t,x)=
\left(\Gamma(h,F_t(x))\right)^{-1}\Gamma(h\circ F_t,x) \\ =
\left(\Gamma(h,F_t(x))\right)^{-1}\Gamma(e^{-At}\circ h,x)   \\ =
\left(\Gamma(h,F_t(x))\right)^{-1}e^{-Ct}\Gamma(h,x) ,
\end{eqnarray*}
we conclude that if $B$ commutes with $C$ and with $\Gamma(h,x)$
for all $x\in\mathbb{D}_1$, then $B$ commutes with all operators
$\widehat{\Gamma}(F_t,x),\ t\ge0,\ x\in\mathbb{D}_1$. Thus, we can
apply Theorem~\ref{th-spiral-gen} with $\Psi_t=e^{-At}$ and
$G_s=e^{-Bs}$ (see Remark~\ref{rem4}). According to this theorem,
for each point $(z,w)\in\Phi[h](\mathbb{D}) $ the image
$\Phi[h](\mathbb{D})$ contains the curve
\[
\biggl\{\left(e^{-At}z,e^{-Ct}G_t(w) \right),\ t\ge0
\biggr\}=\biggl\{\left(e^{-At}z,e^{-(B+C)t}w \right),\ t\ge0
\biggr\}.
\]
So, by Definition~\ref{def-spiral} (see also Remark~\ref{rem7}),
the mapping $\Phi[h]$ is $\left(
\begin{array}{cc}
  A & 0 \\
  0 & B+C \\
\end{array}
\right) $-spirallike. The proof is complete. \epr

\noindent{\bf Proof of Theorem~\ref{th-conv_dir}.} Let
$h\in\Hol(\mathbb{D}_1,X)$ be a mapping convex in the direction
$\tau$ which satisfies~(\ref{assump-c}). Let $B\in L(Y)$ be an
accretive operator which commutes with $C$. Then the semigroup
$\left\{e^{-Bs}\right\}_{s\ge0}$ consists of proper contractions
with respect to the norm $\|\cdot\|_Y$. As in the proof of
Theorem~\ref{th-spiral}, we conclude that if $B$ commutes with $C$
and with $\Gamma(h,x)$ for all $x\in\mathbb{D}_1$, then $B$
commutes with all operators $\widehat{\Gamma}(F_t,x),\ t\ge0,\
x\in\mathbb{D}_1$. Once again, we can apply
Theorem~\ref{th-spiral-gen} with $\Psi_t(z)=z+t\tau$ and
$G_s(w)=e^{-Bs}w$. This theorem implies that for each point
$(z,w)\in\Phi[h](\mathbb{D}_1) $ the image $\Phi[h](\mathbb{D}_1)$
contains the curve
\[
\biggl\{\left(z+t\tau,e^{-Ct}G_t(w) \right),\ t\ge0
\biggr\}=\biggl\{\left(z+t\tau,e^{-(B+C)t}w \right),\ t\ge0
\biggr\}.
\]
The proof is complete. \epr

\section{Concluding remarks}
\setcounter{equation}{0}

\noindent{\bf 1. Bloch type mappings}

\begin{propo}
Let $\Gamma$ be an appropriate operator. Suppose that a mapping
$h\in\mathcal{K}$ satisfies the following conditions:
\begin{itemize}
\item[(i)]
$\displaystyle\sup_{x\in\mathbb{D}_1}\|h'(x)\|_{L(X)}\left(
1-\|x\|_X^2 \right)<\infty$;

\item[(ii)]
$\displaystyle\sup_{x\in\mathbb{D}_1}\left\|\Gamma(h,x)\right\|_{L(Y)}\left(
1-\|x\|_X^2 \right)<\infty$;

\item[(iii)] $\displaystyle\sup_{x\in\mathbb{D}_1} \left\|
\frac\partial{\partial x}\Gamma(h,x) \right\|_{L(X,L(Y))}
p(\|x\|_X) \left( 1-\|x\|_X^2 \right)<\infty$.
\end{itemize}
Then $\displaystyle\sup_{(x,y)\in\mathbb{D}}\left\|\Phi[h]'(x,y)
\right\|_{L(Z)}\left( 1-\|(x,y)\|^2 \right)<\infty$.
\end{propo}

\pr Differentiating $\Phi[h]$ we get
\[
\Phi[h]'(x,y)\left[(z,w)\right]= \left(h'(x)z,
\frac\partial{\partial x}\Gamma(h,x)[z]y + \Gamma(h,x)w  \right).
\]
The direct estimation leads us to
\begin{eqnarray*}
&& \left\| \Phi[h]'(x,y)\right\|_{L(Z)}
=\sup_{(z,w)\in\mathbb{D}}\left\|
\Phi[h]'(x,y)\left[(z,w)\right]\right\| \\ &&\le
\sup_{(z,w)\in\mathbb{D}}\left(\|h'(x)z\|_X + \left\|
\frac\partial{\partial x}\Gamma(h,x) [z]\right\|_{L(Y)} p(\|x\|_X)
+ \left\|\Gamma(h,x)\right\|_{L(Y)}\|w\|_Y\right) \\ &&\le
\|h'(x)\|_{L(X)} + \left\| \frac\partial{\partial x}\Gamma(h,x)
\right\|_{L(X,L(Y))} p(\|x\|_X) +
\left\|\Gamma(h,x)\right\|_{L(Y)}.
\end{eqnarray*}
Therefore,
\begin{eqnarray*}
&&\sup_{(x,y)\in\mathbb{D}}\left\|\Phi[h]'(x,y)
\right\|_{L(Z)}\left( 1-\|(x,y)\|^2 \right) \\ \le &&
\|h'(x)\|_{L(X)}\left( 1-\|x\|_X^2 \right) +
\left\|\Gamma(h,x)\right\|_{L(Y)}\left( 1-\|x\|_X^2 \right) \\ +
&& \left\| \frac\partial{\partial x}\Gamma(h,x)
\right\|_{L(X,L(Y))} p(\|x\|_X) \left( 1-\|x\|_X^2 \right),
\end{eqnarray*}
and the assertion follows. \epr

\vspace{3mm}

\noindent{\bf 2. Open questions}

a. It seems to be possible to repeat a similar construction for
non-linear operators $\Gamma$. At the same time, we know of no
concrete example of an extension operator of the
form~(\ref{RS-abs}) with non-linear $\Gamma$. The question could
be to find such examples.

b. As a rule, the convexity property is more delicate. For
instance, quoting \cite{G-K-book}, we note that it seems to be
difficult to perturb either the extension operator or the domain
without losing the convexity-preserving property. The original
Roper--Suffridge operator (\ref{RS}) preserves the convexity of
the image of the $p$-ball only if $p=2$, i.e., of the Euclidean
ball. On the other hand, if $f$ is convex, then the extended
mapping defined by formula (\ref{G-K-K}) is convex if and only if
$\beta=\frac12$ (see, \cite{G-H-K-S2002}).  So, it is natural to
examine which conditions on $\Gamma$ allow the extension operator
(\ref{RS-abs}) to preserve the convexity.

c. Our scheme does not cover the extension operators introduced by
Muir \cite{M1, M2}. We ask: how to expand it to include his
operators.

\vspace{2mm}

{\bf Acknowledgments.} This research is part of the European
Science Foundation Networking Programme HCAA. The author is very
grateful to Professor D.~Shoikhet and Dr.~M.~Levenshtein for very
helpful remarks.


\begin{thebibliography}{999}

\bibitem{A-E-S} D. Aharonov, M. Elin and D. Shoikhet,
Spirallike functions with respect to a boundary point, {\it Journ.
Math.  Anal. Appl.} {\bf 280} (2003), 17--29.

\bibitem{B-K} D. M. Burns and S. G. Krantz,
Rigidity of holomorphic mappings and a new Schwarz lemma at the
boundary, {\it J. Amer. Math. Soc.} {\bf 7} (1994), 661--676.

\bibitem{E-K-R-S} M. Elin, D. Khavinson, S. Reich and D. Shoikhet,
Linearization models for parabolic dynamical systems via Abel's
functional equation, {\it Ann. Acad. Sci. Fen.} {\bf 35} (2010),
1--34

\bibitem{E-R-S-2004} M. Elin, S. Reich and D. Shoikhet,
Complex Dynamical Systems and the Geometry of Domains in Banach
Spaces, {\it Dissertationes Math. (Rozprawy Mat.)} {\bf 427}
(2004), 62 pp.

\bibitem{E-S3} M. Elin and D. Shoikhet,
Semigroups with boundary fixed points on the unit Hilbert ball and
spirallike mappings, in: {\it Geometric Function Theory in Several
Complex Variables, 82--117}, World Sci. Publishing, River Edge,
NJ, 2004.

\bibitem{F-L} S. Feng and T. S. Liu,
The generalized Roper--Suffridge operator, {\it Acta Mathematica
Scientia} {\bf 28B} (2008), 63--80.

\bibitem{GS-99} S. Gong,
{\it Convex and  starlike mappings in several complex variables},
Science Press, Beijing--New York \& Kluwer Acad. Publ.,
Dordrecht--Boston--London, 1999.

\bibitem{GAW} A. W. Goodman,
{\it Univalent Functions}, Vols. I, II, Mariner Publ. Co., Tampa,
FL, 1983.

\bibitem{G-H-K-S2002} I. Graham, H. Hamada, G. Kohr and T. J. Suffridge,
Extension operators for locally univalent mappings, {Michigan
Math. J.} {\bf 50} (2002), 37--55.

\bibitem{G-K-2000} I. Graham and G. Kohr,
Univalent mappings associated with the Roper-Suffridge extension
operator, {\it J. Analyse Math.} {\bf 81} (2000), 331--342.


\bibitem{G-K-book} I. Graham and G. Kohr,
{\it Geometric Function Theory in One and Higher Dimensions},
Marcel Dekker Inc., New York--Basel, 2003.

\bibitem{G-K-K-2000} I. Graham, G. Kohr and M. Kohr,
Loewner chains and Roper--Suffridge extension operator, {\it J.
Math. Anal. Appl.} {\bf 247} (2000), 448--465.

\bibitem{G-K-P-2007} I. Graham, G. Kohr and J. A. Pfaltzgraff,
Parametric representation and linear functionals associated with
extension operators for biholomorphic mappings, {\it Rev. Roumaine
Math. Pures Appl.} {\bf 52} (2007), 47--68.

\bibitem{L-L-2005} X. S. Liu and T. S. Liu,
The generalized Roper--Suffridge extension operator on a Reinhardt
domain and the unit ball in a complex Hilbert space, {\it Chinese
Ann. Math. Ser. A} {\bf 26} (2005), 721--730.

\bibitem{M1} J. R. Jr. Muir,
A modification of the Roper--Suffridge extension operator, {\it Comput. Methods Funct. Theory} {\bf 5} (2005), 237--251.

\bibitem{M2} J. R. Jr. Muir,
A class of Loewner chain preserving extension operators, {\it J. Math. Anal. Appl.} {\bf 337} (2008), 862--879.

\bibitem{P-S} J. A. Pfaltzgraff and T. J. Suffridge,
An extension theorem and linear invariant families generated by
starlike maps, {\it Ann. Univ. Mariae Curie-Sk\l odowska, Sect, A}
{\bf 53} (1999), 193--207.

\bibitem{RS-SD1} S. Reich and D. Shoikhet,
Generation theory for semigroups of holomorphic mappings in Banach
spaces, \textit{Abstr. Appl. Anal.} \textbf{1} (1996), 1--44.

\bibitem{RS-SD-book} S. Reich and D. Shoikhet,
{\it Fixed Points, Nonlinear Semigroups and the Geometry of
Domains in Banach Spaces}, World Scientific Publisher, Imperial
College Press, London, 2005.

\bibitem{Ro-Su} K. Roper and T. J. Suffridge,
Convex mappings on the unit ball of $\C^n$, {\it J. Analyse Math.}
{\bf 65} (1995), 333-–347.

\bibitem{STJ-77} T. J. Suffridge,
Starlikeness, convexity and other geometric properties of
holomorphic maps in higher dimensions, {\it Complex Analysis
(Proc. Conf. Univ. Kentucky, Lexington, KY, 1976)}, Lecture Notes
in Math. {\bf 599}, 1977, 146--159.


\end{thebibliography}
\end{document}